\documentclass[a4paper]{amsart}
\usepackage{epsfig}
\usepackage{amsmath}
\usepackage{amsthm}
\usepackage{amssymb}
\usepackage{latexsym,amsfonts}
\usepackage{amscd}
\usepackage[all]{xy}
\parskip\smallskipamount
\numberwithin{equation}{section}
\newcommand{\ck}{{{k}}} 

\newcommand{\mand}{\mbox{ and }}
\newcommand{\mwith}{\mbox{ with }}
\newcommand{\mif}{\mbox{ if }}
\newcommand{\mor}{\mbox{ or }}

\newcommand{\mAut}{\mbox{Aut\,}}

\newcommand{\melse}{\mbox{ else }}

\newcommand{\mforall}{\mbox{ for all }}
\newcommand{\mforsome}{\mbox{ for some }}
\newcommand{\mGl}{\mbox{Gl}}

\newcommand{\mHom}{\mbox{Hom}\,}

\newcommand{\mmax}{\mbox{max\,}}

\newcommand{\mmin}{\mbox{min\,}}

\newcommand{\Aa}{{\Bbb A}}

\newcommand{\Rr}{{\Bbb R}}

\newcommand{\Zz}{{\Bbb Z}}

\newcommand{\cM}{{\cal{M}}}

\newcommand{\cO}{{\cal{O}}}
\newcommand{\cP}{{\cal{P}}}
\newcommand{\cR}{{\cal{R}}}

\newcommand{\cT}{{\cal{T}}}

\newcommand{\frp}{{\mathfrak p}}

\newcommand{\lra}{\longrightarrow}

\newcommand{\lla}{\longleftarrow}

\newtheorem{Prop}{Proposition}[section]
\newtheorem{Thm}[Prop]{Theorem}
\newtheorem{Lemma}[Prop]{Lemma}
\newtheorem{Cor}[Prop]{Corollary}

\let\cal\mathcal

\def\mHom{\operatorname {Hom}}
\def\Ext{\operatorname {Ext}}
\def\mExt{\operatorname {Ext}}
\def\mEnd{\operatorname {End}}

\def\rk{\operatorname {rk}}

\def\End{\operatorname{End}}

\def\NA{\operatorname{NA}}

\newcommand{\udim}{\mathrm{\underline{dim}}}

\begin{document}
\title{On the complement of the dense orbit for a quiver of type $\Aa$}
\author{Karin Baur, Lutz Hille}

\date{06.05.2010}
\begin{abstract}
Let $\Aa_t$ be the directed quiver of type $\Aa$ with $t$
vertices. For each dimension vector $d$ there is a dense orbit in the
corresponding representation space. The principal aim of this note is to use just
rank conditions to define the irreducible components in the complement
of the dense orbit. Then we compare this result with already existing
ones by Knight and Zelevinsky, and by Ringel. Moreover, we compare
with the fan associated to  the quiver $\Aa$ and derive a new formula
for the number of orbits using nilpotent classes.  
In the complement of the dense orbit we determine the irreducible
components and  their codimension. Finally, we consider several
particular examples.  
\end{abstract}
\maketitle
\tableofcontents

\section{Introduction}\label{Sintro}

The principal aim of this note is to describe the complement of the
generic orbit in the representation space of a directed quiver of type
$\Aa_t$ with vertices $\{1,2,\ldots,t \}$ and arrows $\alpha_i: i+1 \lra
i$. For a dimension vector $d = (d_1,\ldots,d_t)$ and a
representation $A = (A_{i,i+1}) = (A_{1,2}, A_{2,3}, \ldots, 
  A_{t-1,t})$ with $A_{i,i+1}: V_{i+1} \lra V_i$ we define     
$$
r_{i,j} :=  \min\{d_l\mid i \leq l \leq j\} \mand A_{(i,j)} :=
A_{i,i+1}A_{i+1,i+2} \ldots A_{j-1,j}.
$$                
Let $Y$ be defined as   the complement in
$$
\cR(Q,d) = \bigoplus_{i = 1}^{t-1} \mHom(\ck^{d_i}, \ck^{d_{i+1}}) \mbox{
  with action of } \mGl(d) = \prod_{i=1}^t \mGl(d_i) 
$$
of the generic orbit
$$
\cO(d) := \{ A = (A_{1,2},\ldots,A_{i,i+1},\ldots,A_{t-1,t}) \in \cR(Q,d)\mid
\rk A_{(i,j)} = 
r_{i,j} \}.
$$
In the complement $Y$ we define closed (not neccesarily irreducible) varieties
$$
Y_{i,j} := \{(A_{i,i+1})_{i=1}^{t-1} \in \cR(Q,d) \mid \rk
A_{(i,j)} \leq r_{i,j} - 1\}.
$$ 
We claim in our main result that all irreducible components of $Y$ are
among the $Y_{i,j}$ and at most $t-1$ of the $Y_{i,j}$ occure as
irreducible components in $Y$. 

For the formulation of the main result we need to define a set of pairs 
$$
J(d) := \{ (i,j) \mid 1 \leq i < j \leq t \mand \mforall i < l < j,  d_l
> \max \{ d_i, d_j \} \}. 
$$
and a subset
$I(d)$ consisting of alle elements in $J(d)$ satisfying one of the
following properties (where we define $d_0 = d_{t+1} = 0$) \\
(i) $d_i = d_j$; \\
(ii) $d_i < d_j$ and we define $a$ to be the minimal index $a > j$ with
$d_a < d_i$. Then $d_l \geq d_j$ for all $j < l < a$.
\\
(iii) $d_i > d_j$ and we define $b$ to be the maximal index $b < i$
with $d_b < d_j$. Then $d_l \geq d_i$ for all $b < l < i$.


If $(i,j) \in J(d)$ then we show that $Y_{i,j}$ is irreducible (Theorem
\ref{Tmain}) and there exists a unique representation $M(i,j)$ whose
orbit is dense in  $Y_{i,j}$. For any irreducible component of $Y$
there exists a representation $M$ whose orbit is dense in this
component. Such a representation $M$ is called {\sl almost generic}. 


Then we prove the following theorem.

\begin{Thm}\label{Tmain}
Assume $d_i > 0$ for each $1=1,\ldots,t$. \\
(i) 
$$Y = \cup_{(i,j) \in I(d)} Y_{i,j} $$
is the decomposition of $Y$ into pairwise different irreducible components. 
\\ 
(ii) Each component is the closure of an orbit  corresponding to an
almost generic representation $M(i,j)$. \\
(iii) For any $(i,j) \in J(d)$ the codimension of $Y_{i,j}$ is $|d_j - d_i| + 1$. \\
(iv) The irreducible components $Y_{i,j}$ in $Y$ of codimension $1$
are in 
bijection with the pairs $(i,j)$ with $d_i = d_j$ and $d_l > d_i$
for all $i < l < j$.
\end{Thm}

In fact we prove the following stronger results. First of all
$Y_{i,j}$ is irreducible precisely when $(i,j)$ is in $J(d)$ (Prop
\ref{Pirred}, Cor. \ref{Cirred} ). Then
$Y$ obviously decomposes into the union of all possible $Y_{i,j}$
(Lemma \ref{Ldecomp}). Next we show that any $Y_{k,l}$ for $(k,l)
\notin I(d)$ is already contained in a union of some other $Y_{(i,j)}$
(Prop \ref{Pcontain}). Moreover, we 
interprete our result in terms of multisegments and nilpotent
classes in Section \ref{Sfurther}.

Note that the techniques are similar to the ones in \cite{BH1}, our
case corresponds to $\frp_u(d)/\frp_u(d)'$ therein, 
however, the index sets are different, no case follows from the
other. With some technical modifications on the index sets
one can also handle the case $\frp_u(d)/\frp_u(d)^{(l)}$ for the
remaining values of $l$ in a similar way (Section \ref{Sparabol}).
\medskip

The paper is organized as follows. In Section \ref{Sorbits} we only collect the
details we need for the proof of the main result in Section \ref{Sproof}. Then
we proceed in Section \ref{Sfurther} with some further descriptions related to
tilting modules and trees, the structure of the fan associated to
tilting modules and other 
combinatorial descriptions. The associated simplicial complex of the
fan coincides with the simplicial complex considered by Riedtman and
Schofield (\cite{RiedtmannSchofield2}).  Then, in Section
\ref{Sexample} we consider 
several examples that are 
of interest: convex and concave dimension vectors, pure and generic dimension
vectors, and symmetric ones. In the last section we compare with the
results in \cite{BH1} and mention some generalizations without proofs.
\medskip

We always work over an infinite field $k$, the results here do not depend
on the ground field. For finite fields, one needs to modify the
definition of a dense orbit slightly: an orbit is dense, if it is
dense over the algebraic closure. For a partition $\lambda =
(\lambda_1,\ldots,\lambda_n)$ we denote 
by $C(\lambda)$ the corresponding nilpotent class defined by
$$
C(\lambda) = \{ A \in \mEnd(V) \mid \dim A^l(V) = \mmax\{ j \mid
\lambda_j \geq l \} \}.
$$
All varieties are considered over the
algebraic closure and might be reducible. Also the action of the group
should be understood over the algebraic closure. We will always
identify isomorphism classes of representations of $\Aa_t$ (with
directed orientation) with so-called multisegments defined below. With
$\sharp [i,j]$ we denote the number $j-i+1$ of integers in the interval
$[i,j]$.

{\sc Acknowledgment:}  This work started during a stay of both authors in Oberwolfach. We are
indebted to the Institut for the perfect working conditions. The second author was supported by the DFG
priority program SPP 1388 representation theory.

\section{Description of the Orbits} \label{Sorbits}

In this section we recall some of the various descriptions of the isomorphism
classes of representations of $\Aa_t$ with the directed
orientation that we need in the proof. Moreover, we recall some
well-known facts from the 
classification of tilting modules and compute the extension
groups. We proceed with these descriptions in Section
\ref{Sfurther}. Further related results can be found in
\cite{KnightZelevinsky} and in the classical papers
\cite{AbeasisdelFra} and \cite{AbeasisdelFraKraft}. 

\subsection{Multisegements}

A multisegment $M$ consists of a union of
intervalls $[i,j]$ with $1 \leq i \leq j \leq t$, written as $M = \oplus
[i,j]^{a_{i,j}}$ (since a multisegments represents an isomorphism
class of representations we write 'direct sum' instead of
'union'). The dimension vector of such a multisegemt is 
defined as
$$
\udim M = (d(M)_1,\ldots,d(M)_t); \quad d(M)_l := \sum_{(i,j) \mid i \leq
  l \leq j} a_{i,j}. 
$$

There are natural bijections between the multisegments of dimensions
vector $d$, the isomorphism classes of representations of $\Aa_t$,
and the orbits of the $\mGl(d)$--action on $\cR(Q,d)$. Moreover, for
any dimension vector $d$ there exists a unique multisegment $M(d)$
corresponding to the dense orbit. This multisegment can be constructed
recursively as follows: Define $a_{1,t}$ to be the minimum of the
entries $d_i$ in $d$. Then we consider $d^1 := d -
a_{1,t}(1,\ldots,1)$ and consider the longest interval $[i,j]$ in
$d^1$ with minimal $i$. Then 
$$
d^2 := d^1 - a_{i,j} \udim [i,j] = d^1 -  a_{i,j}
(0,\ldots,0,1,\ldots,1,0,\ldots0)
$$ 
is nonnegative for some 
maximal $a_{i,j}$ and we proceed with $d^2$ instead of $d^1$ in the same
way. Eventually, we obtain a multisegment $M(d)$ with at most $t$
different direct summands. A second way to obtain this multisegment is
described in in \cite{Hvol}, Section 8 (this is a similar, but not the
same, construction as in \cite{BHRR}) as follows: consider the unique
diagram with $d_i$ 
vertices in the $i$th column and connect each vertex in the $i$th column
and the $k$th row with the vertex in the $(i+1)$th column and the
$k$th row (if it exists). Roughly one connects all neighboured
vertices in the same row. The connected components of this diagram are
the direct summands and this diagram represents the multisegment
$M(d)$ (see Section \ref{Sexample} for examples).

{\sc Definition. }
A dimension vector $d$ is {\sl generic} if $M(d)$ contains precisely
$t$ pairwise different segments. A dimension vector is {\sl pure} if
$d_1 = d_t$, $d_l \geq d_1 = d_t$ for each $1 \leq l \leq t$ and this
condition holds recursively for each connected component in the
support of $d(\geq a) := (\max\{d_1-a, 0\}, \ldots, \max\{d_t-a,
0\})$. Examples can be found in Section \ref{Sexample}, see
\ref{Sexample}.1 and \ref{Sexample}.4.  

\subsection{Extensions and homomorphisms}

The category of finite dimensional representations of $\Aa_t$ is a
hereditary category and the Euler characteristic $\langle- ,- \rangle
= \dim \mHom(-,-) - \dim \mExt(-,-)$, respectively the
Hom- and Ext-spaces are 
$$
\begin{array}{ccc}
\langle [i,j], [k,l] \rangle$ is just $\sharp( [i,j] \cap [k,l] ) -
\sharp ([i+1,j+1] \cap [k,l])  \\
\mHom([i,j], [k,l]) = \left\{ 
\begin{array}{ll}
k & \mif k \leq i \leq l \leq j \\
0 & \mbox{ otherwise}
\end{array} \right. \\
\mExt([i,j], [k,l]) = \left\{ 
\begin{array}{ll}
k & \mif i < k \leq j + 1 < l + 1 \\
0 & \mbox{ otherwise}
\end{array}
\right.
\end{array}
$$
All this follows from direct calculations using a projective or
an injective resolution 
$$
0 \lra [j + 1,t] \lra [i,t] \lra [i,j] \qquad \mor \qquad [i,j] \lra
[1,j]  \lra    [1,i+1] \lra 0.
$$
 
\begin{Prop}\label{Pextgroups}
a) A multisegment $M$ has no selfextension precisely when for each
pair of direct summands $[i,j]$ and $[k,l]$ of $M$ one of the
following conditions hold \\
(i) $[i,j] \subseteq [k,l]$, \quad
(ii) $j < k-1$, \quad
(iii) $[k,l] \subseteq [i,j]$, or \quad
(iv) $l < i-1$. \\
b) The multisegment $M(d)$ has no selfextension and any other
multisegment $M$ of dimension vector $d$ satisfies $\mExt(M,M) \not=
0$. \\
c) A multisegement $M$ satisfies $\mExt(M,M) = k$ precisely when it
contains two segments
$[i,j]$ and $[k,l]$ with $j \geq  k-1$, $i < k$, and $j < l$ as a
direct summand and the 
complement of $[i,j]$ equals $M(d')$, where  $d' = \udim M - \udim
[i,j]$ and the complement of $[k,l]$ equals $M(d'')$, where  $d'' =
\udim M - \udim [k,l]$. \\ 
d) A multisegmet $M = \oplus [i,j]^{a_{i,j}}$ is almost generic
precisely when the direct sum of the pairwise non-isomorphic direct
summands $N = \oplus_{(i,j) \mid a_{i,j} > 0} [i,j]^{}$ satisfies
$\mExt^1(N,N) = k$ and one of the direct summands with non-trivial
Ext--group occurs with multiplicity one in $M$.
\end{Prop}

{\sc Proof. }
a) and the first claim of b) is a direct consequence from the formula
for the extension groups above. The uniqueness in b) follows either
directly from the construction, or since $\cR(Q,d)$ is irreducible (it
can contain at most one dense orbit). To prove c) one uses that
$\mExt^1$ is additive, thus there is at most one non--vanishing
extension group. Finally, to prove d) we note that for $d = \udim N$
we have $\dim \mEnd(N,N) =
\dim \mEnd(M(d),M(d)) +1 $ by a simple computation of the Euler characteristic
$\langle M(d),M(d) \rangle = \langle N,N \rangle.$ Thus, the stabilizer of
the orbit of $M(d)$ and $N$ differ by one and then the dimension of the
orbits also differ by one. The closure of the orbit of $M(d)$ obviously
contains the orbit of $N$. Now assume $M$ is a multisegmet as in the
claim and it is neither generic nor almost generic. Take two direct
summands $[a,b]$ and $[c,d]$ of $M$ with $\mExt([a,b],[c,d]) \not=
0$. We define a new multisegment $M'$ of the same dimension vector by
deleting $[a,b]$ and $[c,d]$ and replacing it by $[c,b] \oplus
[a,d]$ (if $c = b+1$ we replace it just by $[a,d]$). Then $M$ is in the closure of the orbit of $M'$ and $M$ is
almost generic precisely when $M'$ is already generic. This in turn is
equivalent to the second condition in d), proving one direction of the
claim.\\
Now assume $N$ satisfies $\mExt(N,N) = k$. Then the closure of the
orbit of $N$ is of codimension one in the space of all representations
of dimensions vector $\udim N$. Thus it is
some irreducible component in the complement of the dense
orbit. Now we add the remaining segments to $N$ so that we obtain $M =
M' \oplus N$. The multisegment $M'$ is then also a direct summand of
$M(d)$, since on can get $M(d)$ from $N \oplus M'$ by extending only
two segments in $N$.

Assume $M$ contains two indecomposable direct summands $[a,b]$
and $[c,d]$, both occuring with multiplicity at least two and
$\mExt([a,b],[c,d]) \not= 0$, then we can again (using extensions)
construct an orbit that is not generic and contains $M$ in ist
closure. Consequently, such an $M$ is not almost generic.
\hfill $\Box$

The proof also follows directly from Zwara's result \cite{Zwara} that the partial
order of the 
Ext-degeneration and the partial order for the geometric degeneration coincide. In the
proof above we only used the trivial direction.

\subsection{Rank conditions}

To any representation $A$ of $Q$ one can associate the ranks of the
compositions of the corresponding matrices. Consider $A = (A_{i,i+1})
\in \cR(Q,d)$. Then we define the {\sl rank triangle}
$$
r(A) = (r_{i,j}(A))_{1 \leq i < j \leq t}, \mwith r_{i,j}(A) = \rk
A_{i,i+1}\cdot \ldots \cdot A_{j-1,j} = \rk A_{(i,j)}.
$$
Moreover, it is convenient to define the {\sl extended rank triangle}
with $r_{i,i} := d_i$ and to define $r_{i,j} = 0$ whenever $i \leq
0$ or $j > t$. 
Obviously, we must have $r_{i,j}(A) \leq r_{i,j} := \mmin\{d_l \mid i
\leq l \leq 
j\}$ and (using generic matrices) the set
$$
\cO(d) := \{ A \in \cR(Q,d) \mid r_{i,j}(A) = r_{i,j} \}
$$
is open and dense in $\cR(Q,d)$. In fact, the set $\cO(d)$ consists of
all representations isomorphic to $M(d)$, since $r_{i,j}(M(d)) =
r_{i,j}$ by construction.

We fix a dimension vector $d$ and 
consider any triangle $s = (s_{i,j})$ of non-negative
integers $s_{i,j}$ satisfying $s_{i,j} \leq r_{i,j}$. Then
$$
X_s^0 := \{ A \in \cR(Q,d) \mid r_{i,j}(A) = s_{i,j} \} \subseteq X_s
:= \{ A \in \cR(Q,d) \mid r_{i,j}(A) \leq s_{i,j} \}  
$$
defines an open (possibly empty) subvariety $X_s^0$ in a closed,
non-empty algebraic subvariety $X_s$ (not necessarily
irreducible) of $\cR(Q,d)$. The rank triangles are partially ordered
by $s \leq u$ iff $u - s$ has only non-negative entries. It turns out
that some of the $X_s$ are irreducible (we determine which ones) and
the rank conditions are very useful for determining the components in
the orbit closures. Moreover, one can reconstruct the multisegment $M$
from the rank condition $s$, where the orbit of $M$ is dense in $X_s$
with $s$ minimal: A direct sum $[i,j]^{a}$ is a
direct summand of $M$ (with maximal possible $a$) if and only if $a =
r_{i,j} - r_{i+1,j} - r_{i,j+1} + r_{i+1,j+1}$. Consequently, $X_s^0$
is empty, if some $r_{i,j} - r_{i+1,j} - r_{i,j+1} + r_{i+1,j+1}$ is
negative. Otherwise $X_s^0$ is dense in $X_s$.

Conversely, given a multisegment $M$ we can easily determine its rank
vector $r(M) 
=(r(M)_{i,j})$ as follows
$$
r(M)_{i,j} = \sharp \left\{ [k,l] \in M \mid k \leq i \leq j \leq l
\right\}.
$$
In the particular case of a segment $[k,l]$, we obtain just the
characteristic function of a triangle as the rank triangle
$$
r([k,l])_{i,j} = 
\left\{
\begin{array}{l}
1 \mif [i,j] \subseteq [k,l] \\
0 \melse.
\end{array}
\right.
$$

\begin{Prop}
a)  If $s \leq u$ then $X_s \subseteq X_u$. In particular, $X_r$
  contains each $X_s$ and $X_0$ (consisting of the zero matrix) is
  contained in each $X_u$. \\
b) The variety $X_s$ is irreducible precisely when it is the closure
  of one $\mGl(d)$--orbit. \\
c) $X_s^0$ is non-empty precisely when $s$ is a sum of functions of
  the form $r([i,j])$ and this is equivalent to $s_{i,j} - s_{i+1,j}
  - s_{i,j+1} + s_{i+1,j+1} \geq 0$ for all pairs $(i,j)$.
\end{Prop}

{\sc Proof. }
Assertion a) is obvious, since $\rk A_{(i,j)} \leq a$ implies $\rk
A_{(i,j)} \leq b$ for any $b > a$.

To prove b) we decompose $X_s$ in a disjoint union of
$\mGl(d)$--orbits. This is possible, since  $X_s$ is
$\mGl(d)$--invariant. Thus we obtain a set of multisegments $\cM_s$ with
$$
X_s = \bigsqcup_{M \in \cM_s} \mGl(d) M.
$$
Consequently, $X_s$ is the union of a finite number of orbit closures
$\overline{\mGl(d) M}$ for a finite number of multisegments $M$. We
can assume this set is minimal. Thus $X_s$ is irreducible precisely
when $X_s = \overline{\mGl(d) M}$ for some maximal $M$ in $\cM_s$.

For c), note that $X_s^0$ is nonempty, precisley when there exists a
multisegment $M$ with $s = \rk M$. This is also equivalent to $s =
\sum_{(i,j)} a_{(i,j)} \rk [i,j]$ is the sum of rank functions of
segments.  To prove the last characterization we note that for $s = 
\sum_{(i,j)} a_{(i,j)} \rk [i,j]$ we obtain $s_{i,j} - s_{i+1,j}
  - s_{i,j+1} + s_{i+1,j+1} = a_{(i,j)} \geq 0$. Conversely, if
  $s_{i,j} - s_{i+1,j}   - s_{i,j+1} + s_{i+1,j+1} \geq 0$ then we
  define $ a_{(i,j)} = s_{i,j} - s_{i+1,j} - s_{i,j+1} +
  s_{i+1,j+1}$. 
 \hfill $\Box$

\section{Proof of the main theorem}\label{Sproof}

We start this section by showing that some of the $Y_{i,j}$ are
irreducible and compute their dimension. Then we show that all
$Y_{i,j}$ for $(i,j)$ not in $I(d)$ are already contained in some union
of other ones. This allows a reduction to the case $Y_{i,j}$ for
$(i,j) \in I(d)$. Finally we show that $Y$ is already contained in the
union of all $Y_{i,j}$.

\subsection{Irreducible varieties}

\begin{Prop}\label{Pirred}
Assume $(i,j) \in J(d)$, then $Y_{i,j}$ is irreducible of
codimension $|d_j - d_i| + 1$ in $\cR(Q,d)$.
 \end{Prop}

{\sc Proof. }
We consider the projection of a representation of $Q$ to the quiver
$Q'$ with vertices $i,i+1,\ldots, j-1,j$ and its subvarieties
$Y_{i,j}$ in $\cR(Q,d)$ and $Y'_{i,j}$ in  $\cR(Q',d')$ defined by
$\rk A_{(i,j)} < r_{i,j}$. Then $\cR(Q,d)$ is a direct product of
$\cR(Q',d')$ with some affine space and $Y_{i,j}$ is a product of
$Y'_{i,j}$ with some affine space. Thus  $Y_{i,j}$ is irreducible
precisely when $Y'_{i,j}$ is irreducible. Consequently, it is
sufficient to prove the claim for $Y_{1,t}$ in  $\cR(Q,d)$. 

We now assume $(i,j) = (1,t)$ and $d_i > d_1, d_t$ for any $1 < i < t$.

Now we consider a multisegment $M$ consisting of $[1,t-1] \oplus
[2,t]$ and $M(e)$ for $e = d - (1,2,2,\ldots,2,2,1)$. A computation of
the ranks $r_{i,j}(M)$ yields  $r_{1,t}(M) = r_{1,t}-1 $ and
$r_{i,j}(M) = r_{i,j} $ for all $(i,j) \not= (1,t)$. Thus the equation
$s_{1,t} = r_{1,t}-1$ and $s_{i,j} = r_{i,j}$ for $(i,j) \not= (1,t)$
defines an orbit $X_s$ and
$X_s$ is the closure of this orbit containing $M$. Consequently it is
irreducible, and it coincides with $Y_{1,t}$.   

Finally, we need to compute the codimension of the orbit closure
$Y_{1,t}$. For this we compute the dimension of the stabilizer of $M(d)$
and of $M$ constructed above. To make the computation easier, we
delete the common direct summands that contribute with the same
dimension to the stabilizer and assume without loss of generality $d_1
\geq  d_t$. Then we need to compute
$$
\begin{array}{c}
\dim \mEnd([2,t-1] \oplus [1,t]^a \oplus [1,t-1]^b ) = a^2 + b^2 +
ab + b + 1\\
\dim \mEnd([2,t] \oplus [1,t]^{a-1} \oplus [1,t-1]^{b+1}) = a^2 +
b^2 + ab + 2b + 2.
\end{array}
$$

Back to $M(d)$, we decompose it into $M(d) =
[2,t-1] \oplus [1,t]^a \oplus [1,t-1]^b \oplus M'$ with maximal $a$ and
$b$. Then $M$ is $[2,t] \oplus [1,t]^{a-1} \oplus [1,t-1]^{b+1} \oplus
M'$ and 
$$
\begin{array}{l}
-(b + 1) = \dim \mEnd(M(d)) - \dim \mEnd(M) = \\ \dim \mEnd([2,t-1] \oplus [1,t]^a
\oplus [1,t-1]^b) - \dim \mEnd([2,t] \oplus [1,t]^{a-1} \oplus
       [1,t-1]^{b+1}). 
\end{array}
$$
Consequently, the codimension of the orbit of $M$ equals $b+1 = d_1
-d_t +1$ and
this equals the codimension of $Y_{1,t}$. Finally, note that under the
reduction from arbitrary $Y_{i,j}$ to $Y_{1,t}$ the codimension does
not change.
\hfill $\Box$ 

\subsection{The reduction process}

\begin{Prop} \label{Pcontain}
a) Assume $(i,j) \notin J(d)$ then there exists some $l$ with $i < l
< j$ and $d_l \leq d_k$ for all $i <  k < j$. In particular,  $d_l \leq
\mmax\{d_i,d_j\}$. In this case we have an inclusion 
$Y_{i,j} \subseteq Y_{i,l} \cup Y_{l,j}$.  \\
b) If $(i,j) \in J(d) \setminus I(d)$ with $d_i \leq d_j$ then
there exists an $l$ with $l > j$ and  $d_i \leq d_l < d_j$. In this
case we obtain $Y_{i,j} \subset Y_{i,l}$. \\
c) If $(i,j) \in J(d) \setminus I(d)$ with $d_i \geq d_j$ then
there exists an $l$ with $l < i$ and  $d_j \leq d_l < d_i$. In this
case we obtain $Y_{i,j} \subset Y_{l,j}$.
\end{Prop}

{\sc Proof. }
Without loss of generality we may assume $d_i \leq d_j$ in the proof. \\ 
a) Consider the maps $A_{(i,j)}: V_i \lra V_j$, $A_{(i,l)}: V_i \lra
V_l$, and  $A_{(l,j)}: V_l \lra V_j$. We consider two cases.

$d_i \geq d_l$: \\
Assume $\rk A_{(i,j)} < d_l$,
then $\rk A_{(i,l)} < d_l$ or $\rk A_{(l,j)} < d_l$ since $A_{(i,j)}:
    V_i \lra V_l \lra V_j$ factors through $V_l$ with $\dim V_l \leq \dim
    V_i, \dim V_j$.

$d_i < d_l$: \\
Assume $\rk A_{(i,j)} < d_i$, then  $\rk A_{(i,l)} < d_i$

b) Consider the maps $A_{(i,j)}: V_i \lra V_j$ and $A_{(i,l)}: V_i
\lra V_l$. Since $(i,j) \in J(d) \setminus I(d)$ there exist some
$l>j$ with $d_i \leq d_l < d_j$ and $(i,l) \in J(d)$. Then from $\rk
A_{(i,j)} < d_i$ follows $\rk A_{(i,l)} < d_i$.

c) This case is opposite to case b). \hfill $\Box$  


\begin{Lemma}\label{Ldecomp}
$$Y = \bigcup_{1 \leq i < j \leq t} Y_{i,j} = \bigcup_{(i,j) \in
  J(d)} Y_{i,j} = \bigcup_{(i,j) \in
  I(d)} Y_{i,j}.
$$
\end{Lemma}

{\sc Proof. }
The dense orbit is defined by the condition $\rk A_{(i,j)} =
r_{i,j}$. Thus, the complement satisfies $\rk A_{(i,j)} < r_{i,j}$ for
at least one pair $(i,j)$ with $r_{i,j} > 0$. Since $r_{i,j} =
\mmin\{d_i\mid 1 \leq i \leq t \} > 0$ we finish the proof of the
first equality. 

To prove the second one we use the proposition above. From Proposition
\ref{Pcontain} a) we obtain $\bigcup_{1 \leq i < j \leq t} Y_{i,j}
\subseteq \bigcup_{(i,j) \in J(d)}Y_{i,j}$ and from part b) and c) $
\bigcup_{(i,j) \in J(d)}Y_{i,j} \subseteq  \bigcup_{(i,j) \in
  I(d)}Y_{i,j}$. 
\hfill $\Box$

\begin{Cor}\label{Cirred}
The variety $Y_{i,j}$ is irreducible precisely when $(i,j) \in J(d)$.
\end{Cor}

{\sc Proof. }
Thanks to Proposition \ref{Pirred}, we only need to prove that
$Y_{i,j}$ is not irreducible for $(i,j)$ not in $J(d)$. Take $(i,j)$
not in $J(d)$, thus there 
exists an $l$ with $i < l < j$ and $d_l \leq \mmax\{d_i, d_j\}$ is
minimal. Then
by Proposition \ref{Pcontain} a) we have $Y_{i,j} \subseteq Y_{i,l}
\cup Y_{l,j}$. Assume first that $d_l \leq \min\{d_i,d_j\}.$ We claim that $Y_{i,j} = Y_{i,l}
\cup Y_{l,j}$ and this is a proper decomposition, none contains the
other. To see the equality, consider any element $A$ in $ Y_{i,l}
\cup Y_{l,j}$. Then  $\rk(A)_{i,l} < r_{i,l}$ or $\rk(A)_{l,j} <
r_{l,j}$. From each of the inequalities follows $\rk(A)_{i,j} <
r_{i,j}$. On the other hand, there exists a representation with
$\rk(A)_{i,l} = r_{i,l}$ and $\rk(A)_{l,j} <
r_{l,j}$ and vice versa, proving also the last claim. \\
In the second case $ \max{d_i,d_j} \geq d_l > \min\{d_i,d_j\}.$ We
construct two different subvarieties that contain $Y_{i,j}$ and none
contains the other. To
simplify the arguments, we assume without loss of generality $d_i \geq
d_l > d_j$ and, using the first case, $d_l$ is the minimal entry of
$d$ between $d_i$ and $d_j$.  The
first variety is just the orbit closure $Y_{l,j}$, the second one is defined by
$r_{i,l}(A) < r_{i,l}$ and   $r_{i,j}(A) < r_{i,j}$. Using
multisegments (or rank conditions) one can show that we obtain at
least two irreducible components in this way ($Y_{l,j}$ is
irreducible, the other variety need not to be). Anyway, we obtain at
least two irreducible components.
\hfill $\Box$

\section{Further descriptions}\label{Sfurther}

In this section we proceed with the various descriptions of the
irreducible components and the tilting modules started in Section
\ref{Sorbits}.  In particular, we use trees and fans to describe the
irreducible components
 and we relate our description to the nilpotent
class representations defined in \cite{Hhabil}.

\subsection{Trees and tilting modules}\label{Strees}

Let $T$ be a $3$--regular tree with one root and  $t+1$  leaves, where
the leaves are enumerated by $0,1,2,\ldots,t-1,t$. We
denote the set of those trees with $\cT_t$. With $\cT^1_t$ we denote all
trees that have precisely on vertex with four neighbours, all other
vertices have three neighbours and admit one root and $t+1$
leaves. There is a natural map from $\cT^1_t$ to the set of unordered
pairs $\cP^2(\cT_t)$ of
 of trees in $\cT_t$ by ``resolving'' the vertex with four
neighbours and replacing it by two $3$--regular vertices (see Figure 1).  
\medskip

\centerline{\includegraphics[height=2cm]{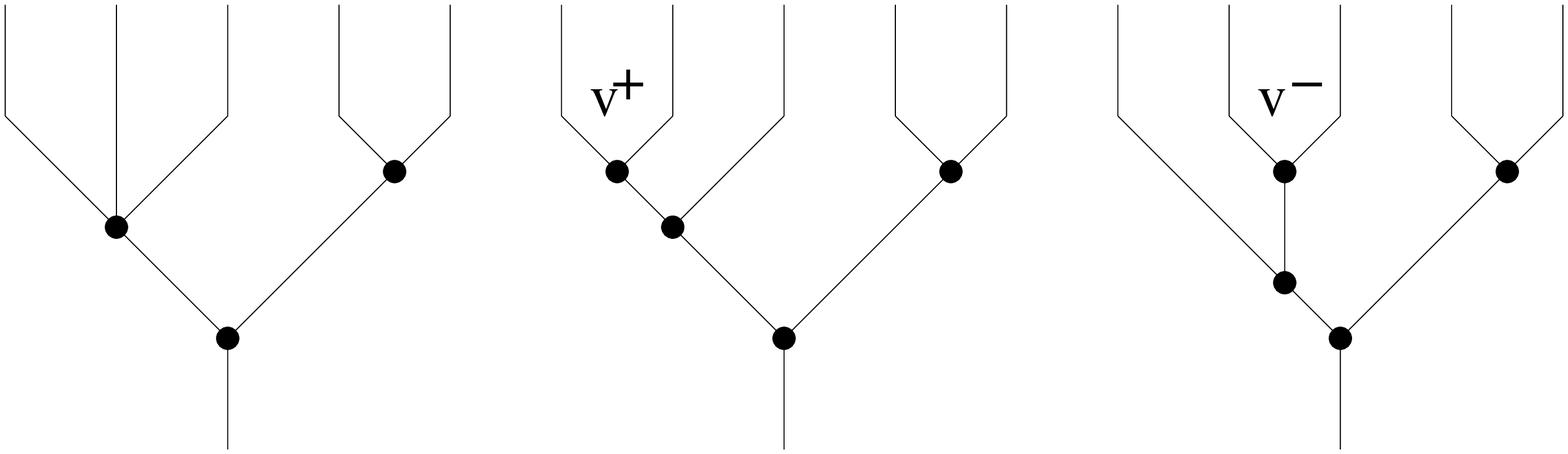}}

\centerline{{\bf Figure 1.} a tree in $\cT^1_4$ and the two associated
  $3$--regular trees in $\cT_4$.}
\medskip

We always draw a tree in the plane and fix the numbering of the leaves
$0,\ldots,t$ from left to right. Two trees are considered to be equal,
if the abstract graphs are isomorphic and the numbering of the leaves
is preserved under the isomorphism. Then each vertex $v$ defines the set of
leaves (in fact an interval) above the vertex $\{ i_T(v)-1, i_T(v),
\ldots,j_T(v) \}$. This way, each vertex $v$ defines a segment $[i(v),j(v)]$.
To any tree in $\cT$ or $\cT^1$ we can associate multisegments as
follows. Assume $T \in \cT_t$ and denote by $T_0$ the vertices in $T$, then we define 
$$
M_T = \bigoplus_{v \in T_0} [i_T(v),j_T(v)]
$$ 
to be the union of the multisegments $[i_T(v),j_T(v)]$
above of $v$. If $S \in \cT^1$ and $T^+$ and $T^-$ are the two
associated $3$--regular trees with unique vertex $v^+ \in T^+$ and
$v^- \in T^-$ (these are the only vertices defining a segment that is
not obtained from the other tree), we define 
$$
\begin{array}{rl}
M_S & = \bigoplus_{v \in  T^+_0}([i_{T^+}(v),j_{T^+}(v)]) \oplus
[i_{T^-}(v^-),j_{T^-}(v^-)] \\
& =  \bigoplus_{v \in T^-_0}([i_{T^-}(v),j_{T^-}(v)]) \oplus
[i_{T^+}(v^+),j_{T^+}(v^+)], \mand \\
\underline M_S & = \bigoplus_{v \in  T_0}[i_{T}(v),j_{T}(v)]
\end{array}
$$
The module $M_T$ for $T \in \cT^1$ has $t+1$ pairwise nonisomorphic
direct summands and the module $\underline M_S$ has $t-1$ pairwise
nonisomorphic direct summands. 

\begin{Thm}
a) If $T$ is a $3$--regular tree, then $M_T = M(d(T))$ for $d(T) =
\udim M_T$. In 
particular, $\Ext(M_T,M_T) = 0 $ for any tree $T$ in $\cT_t$. \\
b) If $S$ is in $\cT^1$ then $\mExt(M_S,M_S) = k$. In particular,
$M_S$ is almost generic. \\
c) If $d$ is generic, then there exists a unique $T \in \cT_t$ so that
$M(d)$ and $M_T$ have the same indecomposable direct summands. Thus
$M(d)$ is a direct summand of several copies of $M_T$.\\
d) If $M$ defines the open dense subset in an irreducible component
$Y_{i,j}$ of
$Y$, then there exists some tree $S \in \cT^1_t$ so
that $M$ is a direct summand of copies of $M_S$.  \\
e) For each multisegment $M$ with $Ext(M,M) = k $ there exists some
$S \in \cT^1$ with $[i_{T^-}(v^-),j_{T^-}(v^-)]$ and
$[i_{T^+}(v^+),j_{T^+}(v^+)]$ as direct summand of $M$.
\end{Thm}   

{\sc Proof. }
Using Prop. \ref{Pextgroups} a) one sees immediately that the
segments in $M_T$ satisfy the vanishing condition for the extension
groups. Thus, the only non-vanishing extension group in $M_S$ are in
the complement of $\underline M_S$, that consists of two
segments. This proves a) and b) (see also the proof in
\cite{Hvol}). 

Part c) also follows from the arguments in loc.~cit.:
Each multisegment with non-vanishing extension group can be completed
to one with $t$ non-isomorphic direct summands. Finally, any
multisegment with precisely $t$ indecomposable summands, all pairwise
non-isomorphic and vanishing extension group is isomorphic to
$M_T$ for some $T$ in $\cT_t$ and the segments determine $T$ uniquely. \\
Using the description of an almost generic multisegment $M$ in
Prop. \ref{Pextgroups} d) we find two segments $[a,b]$ and $[c,d]$
in $M$ with non-vanishing extension group. Moreover, we can assume
that $M$ has $t+1$ non-isomorphic direct summands (otherwise we add
further ones to $M$).  Deleting all summands of the form
$[a,b]$ defines a tree $T^+$, deleting the other direct summand
$[c,d]$ defines a different tree $T^-$ by c). By construction, both
trees come from a common $S$ in $\cT^1_t$ so that $M_S$ and $M$ contain
the same indecomposable direct summands up to isomorphism. This proves
d). The two
summands in e) are just $[a,b]$, respectively $[c,d]$.
\hfill $\Box$

\subsection{Nilpotent class representations}

There is an obvious formula for the number $N(d)$ of orbits in $\cR(Q,d)$. We
just count the number of multisegments $\oplus [i,j]^{a_{i,j}}$
defined by a non-strict triangle $a = (a_{(i,j)})_{1 \leq i\leq j \leq
  t}$ 
$$
N(d) = \sharp \{ a = (a_{i,j})_{1 \leq i \leq j \leq t}\mid a_{i,j}
\in \Zz_{\geq 0}, \quad \sum_{i,j} a_{i,j} \udim [i,j] = d \}.
$$ 
This function is also called {\sl Kostants partition function} for
type $\Aa$.
It is for large $d$ not efficiently computable, thus an easier formula is
desirable. For we define numbers $NA(\lambda, \mu)$ for any two
partitions $\lambda$ of $b > 0$ and $\mu$ of $c > 0$
$$
\NA(\lambda, \mu) = \left\{
\begin{array}{ll}
\prod_{l=1}^{\infty} (\sharp \{i\mid \lambda_i = \mu_i = l \} + 1) &
\mif  |\lambda_i - \mu_i| \leq  1 \mforall i \\
0 & \mif |\lambda_i - \mu_i| \geq 2 \mforsome i.
\end{array}
\right.
$$

\begin{Prop}\label{PNA}
The number of multisegments coincides with the sum, taken over all
sequences of partitions $(\lambda^1, \ldots, \lambda^t)$ with
$\lambda^i$ a partition of $d_i$, $\lambda^1 = (1)^{d_1}$ and
$\lambda^t = (1)^{d_t}$ are both trivial, of the product of the
numbers $\NA(\lambda^i, \lambda^{i+1})  $
$$
N(d) = \sum_{(\lambda^1, \ldots, \lambda^t)} \prod_{i=1}^{t-1}\NA(\lambda^i, \lambda^{i+1}).
$$
\end{Prop}

{\sc Proof. }
We only mention the idea of the proof, the details can be found in
\cite{Hhabil}, Section 4.2. First we consider the preprojective
algebra $\Pi_t$ of $\Aa_t$ and the cyclic quiver $\widetilde \Aa_1$
with two vertices together with the natural projection maps
$$
\times_{i=1}^{t-1} \cR(\widetilde \Aa_1,(d_i,d_{i+1})) \stackrel{p_1}{\lla} \cR(\Pi_t,d)
\stackrel{p_2}{\lra} \times_{i=1}^t \cR(k[T]/T^{d_i}). 
$$
If we denote an element in $\cR(\Pi_t,d)$ by $(A,B)$, then it satisfies
$B_1A_1 = A_1B_1 - B_2A_2 = \ldots =  A_{t-1}B_{t-1} - B_tA_t = A_tB_t
= 0$, where $A_i:V_i \lra V_{i+1}$ and $B_i: V_{i+1} \lra V_i$. The
projections are defined by
$$
p_1:(A,B) \mapsto ((A_1,B_1), \ldots,(A_{t-1},B_{t-1})) \mand $$ 
$$p_2: (A,B) \mapsto
(B_1A_1, A_1B_1, \ldots,A_{t-1}B_{t-1}).
$$
In particular, each element $(A,B)$ defines a sequence of partitions
$(\lambda^1,\ldots, \lambda^t)$ (defined by the partition of the
nilpotent class of $B_1A_1,A_1B_1,\ldots,A_{t-1}B_{t-1}$). By
definition, $\lambda^1$ and $\lambda^t$ are always the trivial ones
corresponding to the zero matrix. It is known (see \cite{KnightZelevinsky} or
\cite{Piasetsky}) that $\cR(\Pi_t,d)$ is equidimensional and the irreducible
components are in bijection with the $\mGl(d)$--orbits on
$\cR(\Aa_t,d)$. Thus $N(d)$ is just the number of irreducible
components in   $\cR(\Pi_,d)$. Now we determine the irreducible
components in a
different way using the projection map above. First note, that
$\NA(\mu, \lambda)$ is the number of irreducible components of
$$
\cR(\widetilde \Aa_1,(\mu, \lambda)) = \{(A,B)\mid AB \in C(\mu), BA
\in C(\lambda)).
$$
If one fixes the sequence of
partitions $(\lambda^1,\ldots,\lambda^t)$, then one can compute the
number of irreducible components in $\cR(\Pi_t,d)$  as the sum of the
products $\NA(\lambda^1, \lambda^2)\cdot \ldots \cdot \NA(\lambda^i, \lambda^{i+1})$ taken over all
such sequences of partitions with $\lambda^1$ and $\lambda^t$ trivial.
\hfill $\Box$

The advantage of the formula above is twofold. First, it is
independent of the orientation of the quiver. We can, for any
orientation of the quiver of type $\Aa_t$ define such a sequence of
partitions. Secondly, the formula in Prop. \ref{PNA} is much more
efficient than just a simple counting. \\
Note that for the generic
representation  $M$ for a a quiver of type $\Aa_t$ with an arbitrary
orientation the corresponding sequence of partitions is just the
trivial one (all $\lambda^i$ are zero). 





\subsection{The fan and the volume}

The sets of trees $\cT_t$ and $\cT^1_t$ define a graph $\Gamma_t$ that is the
dual graph of the simplicial complex of tilting modules defined in
\cite{RiedtmannSchofield2}. This simplicial complex has a natural
realisation as a fan $\Sigma$ in the positive quadrant  $K_{\Rr}^+$ of the real
Grothendieck group $K_0$, where $K_{\Rr}^+: = \Rr_{\geq 0}^t \subset \Rr^t \simeq
K_0 \otimes 
\Rr$. This fan is described in \cite{Hvol}. From the fan, one can
again determine the irreducible compenents in a simple way.
\medskip

We start to define the graph $\Gamma = \Gamma_t$. The vertices
$\Gamma_0$ are just the trees in $\cT_t$. The set of edges is
$\cT^1_t$. The end points of the edge $S$ consists of the two
resolutions $T^+$ and $T^-$ of $S$. 
\medskip

Then we recall the definition of the fan $\Sigma$. For a precise
definition of a fan,
some first properties and applications we refer to
\cite{Fultontoric}. Note first, that a fan
$\Sigma$ is a
finite collection of rational, convex, strongly convex, polyhedral
cones that satisfy two conditions: \\
F1) each face of a cone in $\Sigma$ is in $\Sigma$ and \\
F2) the intersection of two cones in $\Sigma$ is a face of both.
\medskip

Note that we only need finite fans,
for tame and wild quivers one needs to allow also infinite ones.
For each $T \in \cT_t$ we define a cone $\sigma_T \subset K_{\Rr}^+$
as the cone spanned by the dimension vectors of the indecomposable
direct summands of $M_T$
$$
\sigma_T := \{ \sum \Rr_{\geq 0}\udim [i,j] \mid [i,j] \mbox{ is a direct
  summand of } M_T \}.
$$
Two cones $\sigma_{T^+}$ and $\sigma_{T^-}$ have a common facet,
precisely when there exists a tree $S \in \cT^1_t$ with corresponding
trees $S^+ = T^+$ and $S^- = T^-$. The fan $\Sigma$ consists of all
cones, generated by dimension vectors of indecomposable direct
summands of a rigid multisegment $M$ (that is $\mExt(M,M) = 0$).
Already the cones $\sigma_T$ determine the fan
$\Sigma$ consisting of all the cones $\sigma$ that are faces of a cone
$\sigma_T$ (including the cones $\sigma_T$ themself). We recall the
main result from \cite{Hvol} together with some 
easy consequences.

\begin{Thm}\label{Tfan}
a) The cones $\sigma \in \Sigma$ are all generated by a part of a
$\Zz$--basis (they are {\sl smooth} cones). \\
b) The union of the cones $\sigma_T$ (that is the same as the union of
all cones in $\Sigma$) cover $K_{\Rr}^+$. For two cones in $\Sigma$
their intersection is a face of both (and it is also in
$\Sigma$). Each cone is a face of a $t$--dimensional cone and each
$t$--dimensional cone equals $\sigma_T$ for some $T \in \cT_t$. \\
c) A dimension vector $d$ is generic, precisely when it is in the
interior of some cone $\sigma_T$. Consequently, for $d$ generic, $T$
is uniquely determined by $d$. \\
d) For each dimension vector $d$ there exists a unique cone $\sigma
\in \Sigma$ with $d \in \sigma$ and no face of $\sigma$ does contain
$d$. This is equivalent to saying that, $d$ is an element of the relative
interior of the cone $\sigma$. Moreover, the cone $\sigma$ is
generated as a cone by the dimension vectors of the indecomposable
direct summands of $M(d)$. In particular, the dimension of $\sigma$ is
the number of pairwise non-isomorphic indecomposable direct summands
of $M(d)$.\\
e) The dual graph of the $t$--dimensional cones and the
$(t-1)$--dimensional cones, occuring as an intersection of two
$t$--diemensional cones, is $\Gamma_t$. \\
f) Each $t$--dimensional cone has precisely $t-1$ neighbours, that is
$\Gamma_t$ is $(t-1)$--regular.
\end{Thm}

The proof can be found in \cite{Hvol}.

\subsection{Irreducible components and the fan $\Sigma$}

Using the fan, we can again determine the irreducible components in
$Y$. Note that $d$ is contained in some maximal set of cones. We
denote the set of trees $T \in \cT_t$ with $d \in \sigma_T$ by $\cT(d)$. 

Assume $d$ is in the relative interior of a
facet $\sigma_S$ that is the intersection of the two $t$--dimensional
cones $\sigma_{T^+}$ and $\sigma_{T^-}$. Note that $S$ is a tree in
$\cT^1_t$ (compare with section \ref{Strees}). Then $S$ defines two
segments $[a_S,b_S]$ and $[c_S,d_S]$ defined as the unique segments
not in $\sigma_S$ but one in $\sigma_{T^+}$ and the other one in $\sigma_{T^-}$.

Then we obtain one component
just by adding to $[a_S,b_S] \oplus [c_S,d_S]$ the unique generic
complement: $M(d') \oplus [a_S,b_S] \oplus [c_S,d_S].$ In this way, each inner 
facet $\sigma_S$
defines a unique irreducible component $Y_S$.

If $d$ is generic, that is $\cT(d)$ consists of just one tree $T$,
then the components in $Y$ correspond to the $t-1$ neighboured
cones. In fact, each neighboured cone of $\sigma_T$ has a common facet
$\sigma_S$ with 
$\sigma_T$ and the component constructed above defines also a
component for the dimension vector $d$. In this way, we obtain
precisely $t-1$ components. It remains to show that they are pairwise
different. Decompose $M(d) = [1,t]^a \oplus \bigoplus [i,j]^{a_{i,j}}$
into $t$ indecomposable, pairwise nonisomorphic, direct
summands. Since there is a unique cone $\sigma_T$ containing $d$,
$M(d)$ has precisley $t$  pairwise nonisomorphic, direct
summands by Prop.~\ref{Tfan}. Then for each segment $[k,l] \not = [1,t]$
with $a_{k,l} \not= 0$ there exists a unique $S$ so that $[k,l]$
coincides with one of the segments, say $[a_S,b_S]$, constructed
above. Then we obtain a component $Y_S$ just by adding to $[a_S,b_S]
\oplus [c_S,d_S]$ the unique generic complement. Since $[a_S,b_S]$ and
$ [c_S,d_S]$ determine $S$, we obtain the desired $t-1$ irreducible components.


\begin{Thm} \label{ThmGeneric}
Let $d$ be generic with $d$ in the interior of a $t$--dimensional
cone $\sigma \in \Sigma$ , then the irreducible components $Y_{i,j}$ of
$Y$ are in bijection with the set of trees $T$ that define a cone with
a common facet $\sigma_S$ with $\sigma$. For such a tree, take the
two unique segments $[a_S,b_S]$ and $[c_S,d_S]$ (with $a_S < c_S \leq
b_S + 1 < d_S + 1$) in $\sigma_T \cup
\sigma$ with nontrivial extension.  Then $(i,j) = (a_S,d_S)$. \\
 \end{Thm}

\section{Examples} \label{Sexample}
In general, $M(d)$ can have less than $t$ pairwise non-isomorphic
direct summands. Also, the codimension of the components and the number
of components in $Y$ can vary. We discuss several examples, where we
have more precise results. This includes the pure case (all
components have codimension one), the generic case (which contains all
dimension vectors that do not lie on a proper face of a cone in the
fan $\Sigma$),
the concave case, and, eventually,  the convex case.

\subsection{Generic dimension vectors}

 For $d$ generic, we have always $t$
indecomposable direct summands and $d$ lies in the interior of a cone
$\sigma_T$ in the fan $\Sigma$ for some $T \in \cT_t$. 

\begin{Prop}
Assume $d$ is a generic dimension vector. Then $Y$ consists of $t-1$
irreducible components, all have codimension at least $2$.
\end{Prop}

{\sc Proof. } This follows directly from Theorem \ref{ThmGeneric}. 
\hfill $\Box$

\subsection{Concave dimension vectors}

In the concave case (that is $d_1 \geq  d_2 \geq \ldots d_{a-1} \geq
d_a \leq d_{a+1} \leq \ldots \leq d_t$) there are also always $t-1$
components. Moreover $I(d)$ can easily be described.

\begin{Prop} Assume $d$ is a concave dimension vector with $d_i > 0$. \\
a) The sets $I(d)$ and $J(d)$ both coincide with $\{(1,2),(2,3),\ldots
(t-2,t-1),(t-1,t)\}.$\\
b) There are precisely $t-1$ irreducible components $Y_{i,i+1}$ for
$i=1,\ldots, t-1$. The 
codimension of $Y_{i,j}$ equals $|d_i - d_{i+1}| + 1$.\\
c) The variety $Y_{i,j}$ is irreducible precisely when $j = i + 1$. \\
d) The dimension vector $d$ is generic precisely when for all $i < j$
with $d_i = d_j$ there exists some $l$ with $i < l < j$ and $d_l < d_i
= d_j$. \\
e) The components in $Y$ of codimension one correspond to pairs
$(i,i+1)$ with $d_i = d_{i+1}$.
\end{Prop}  

{\sc Proof. } We use our main theorem  together with
the methods obtained in Section 4. \hfill$\Box$

{\sc Example. }
We consider $d = (d_1,\ldots,d_7) =(5,4,3,1,2,4,6)$ and get $I(d) = \{
(1,2),(2,3), \ldots ,(5,6),(6,7) \}$. In this case all components have
codimension at least $2$ and $d$ is generic. The corresponding
multisegment $M(d)$ is $[1,7] \oplus [1,3]^2 \oplus [1,2] \oplus [1,1]
\oplus [5,7] \oplus [6,7]^2 \oplus [7,7]^2.$
\medskip 

\centerline{\includegraphics[height=2cm]{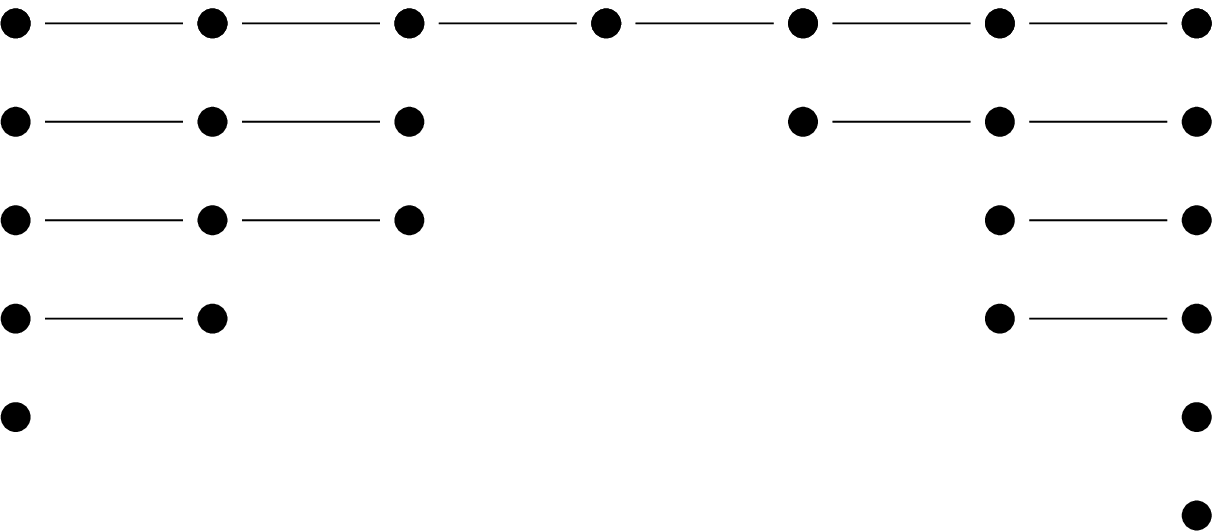}}

\centerline{ {\bf Figure 2.} $M(5,4,3,1,2,4,6)$}
\medskip

and the almost generic multisegments (corresponding to minimal
degenerations) look like (all are different and minimal) 
\medskip

\centerline{\includegraphics[width=10cm]{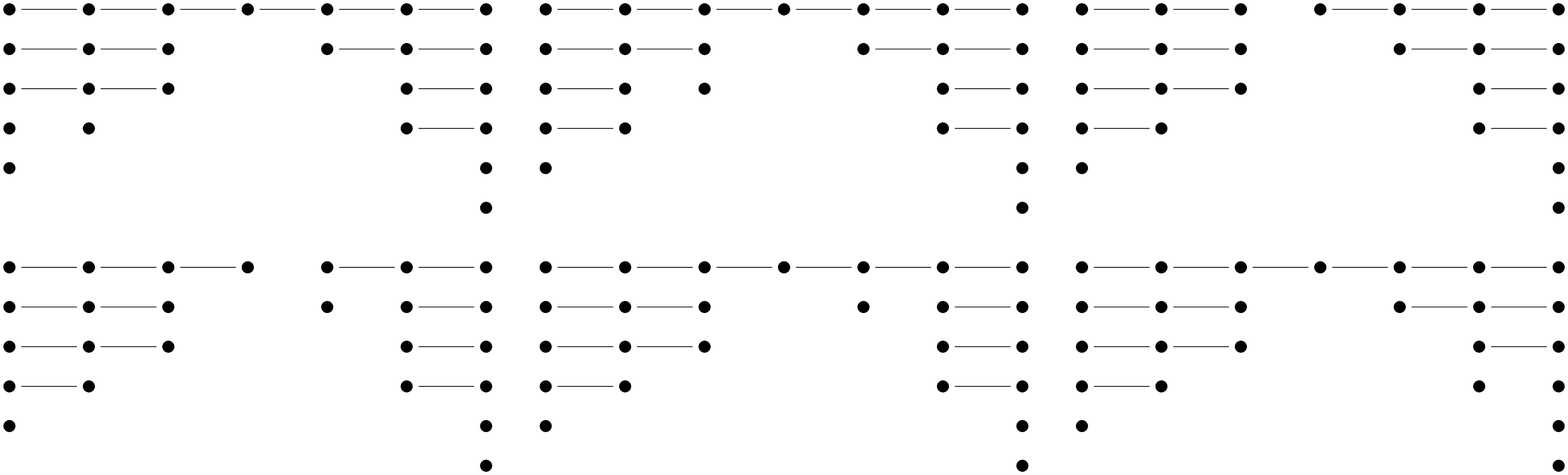}}

\centerline{ {\bf Figure 3.} The minimal degenerations of $M(5,4,3,1,2,4,6)$}
\medskip

Thus we get the following irreducible components

$Y_{i,i+1} := \{(A_{i,i+1}\mid \rk A_{i,i+1} < \mmin \{ d_i,d_{i+1}
\} \}$ of codimension $2$ and $3$, respectively.

\subsection{Unimodular (Convex) Dimension Vectors}

For an unimodular dimension vector (or a convex one) we have the
opposite inequalities $d_1 \leq
d_2 \leq \ldots \leq d_{a-1} \leq d_a \geq d_{a+1} \geq \ldots \geq d_t$  
In the unimodular case there can be less than $t-1$ components. To
illustrate this 
we give two examples, one with $t-1$ components and one with $(t-1)/2$
components, which is the minimal number of components. Moreover, it is
convenient to exclude the previous case of a concave dimension vector
in what follows. We denote by 
$d^+$ the maximal entry in $d$.

\begin{Prop} Let $d$ be an unimodular, sincere dimension vector that
  is not concave.
a) $I(d) \not= J(d)$. \\
b) There are precisely $t-1$ irreducible components in $Y$ if and only
if from $d_i = d_j$ follows $i = j-1$ or $ j+1$.
\end{Prop}

{\sc Proof. }
Under the assumptions we have a maximal entry $d^+$ in $d$ with $d_1
\not= d^+$ and $d_t \not= d^+$. Then $(1,2),(2,3), \ldots,(t-1,t)$ are
all in $J(d)$. If they are also in $I(d)$, then $I(d)$ is just the set
of all pairs $(i,i+1)$. Then we get $d_i = d_{i+1}$
for all $i$, a contradiction to our assumption.
\hfill $\Box$
\medskip

{\sc Example 1.}
We consider $d = (d_1,\ldots,d_7) =(1,3,5,7,6,4,2)$. The corresponding
multisegment is  
\medskip

\centerline{\includegraphics[height=2cm]{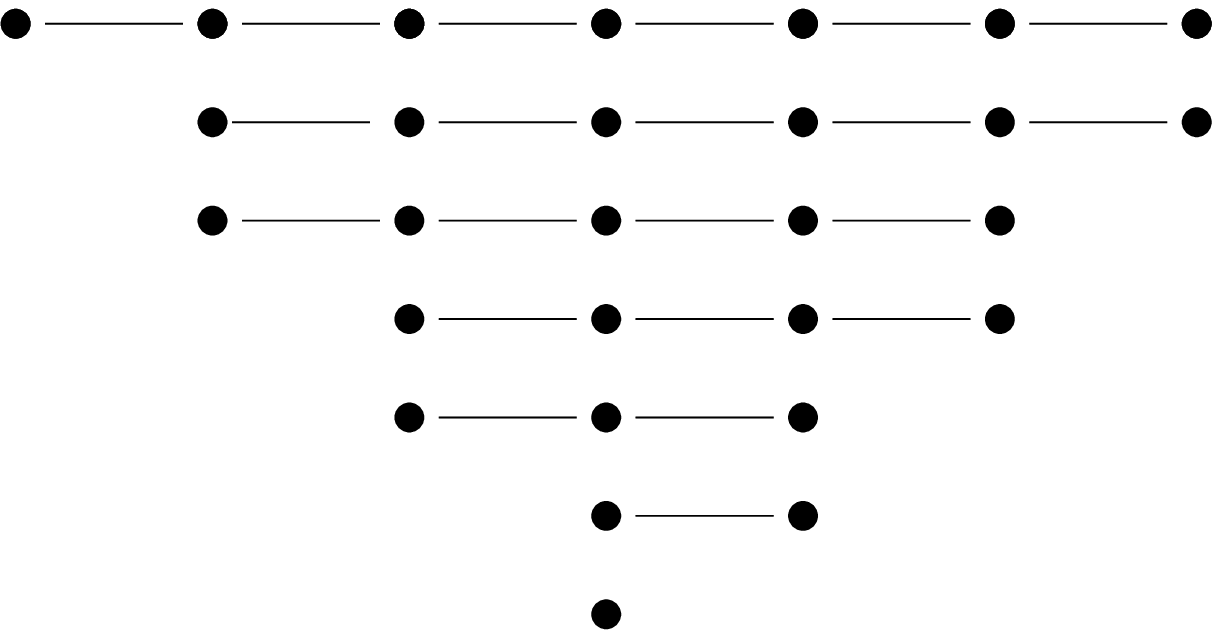}}

\centerline{ {\bf Figure 4.} $M(1,3,5,7,6,4,2)$}
\medskip

and its minimal degenerations look like (all are different and irreducible)
\medskip

\centerline{\includegraphics[width=10cm]{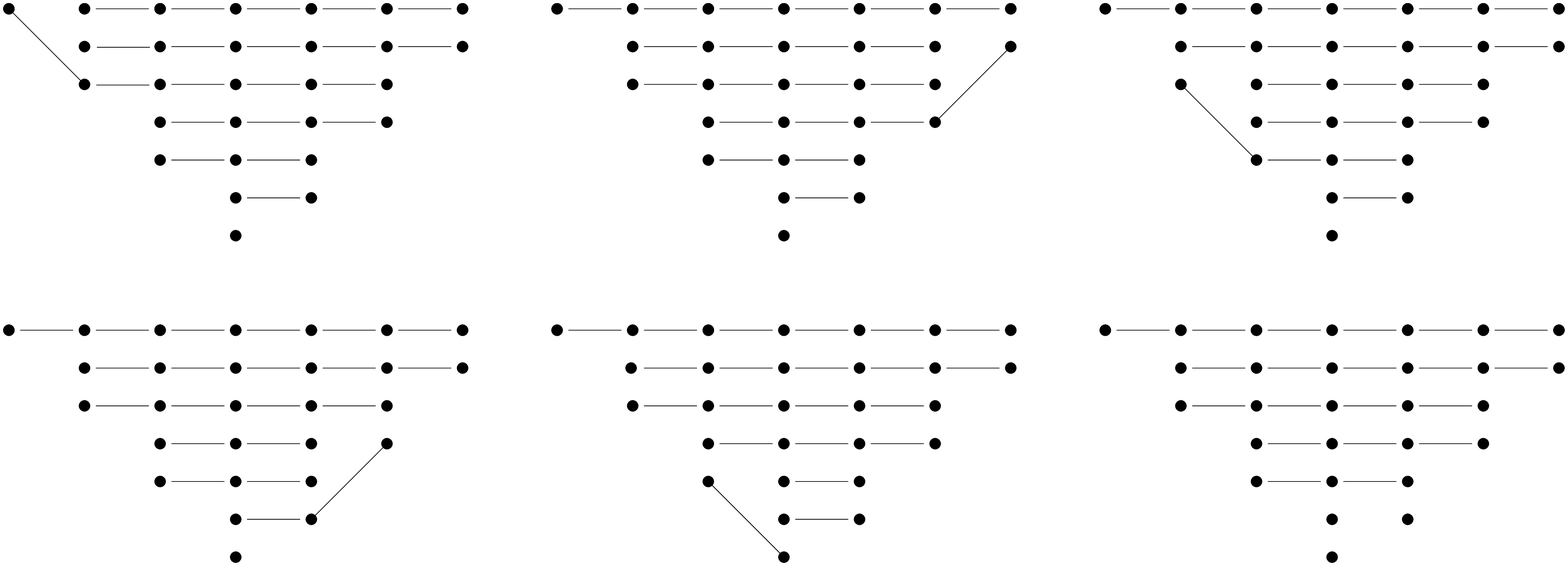}}

\centerline{ {\bf Figure 5.} minimal degenerations of $M(1,3,5,7,6,4,2)$}
\medskip

{\sc Example 2.}
The other extreme case is to have only $(t-1)/2$ components. This
forces that $d$ has often the same entry and the entries in between
two equal entries are all larger. Thus, in the convex case we may consider
$d = (1,2,4,5,4,2,1)$ with $I(d)=\{ (1,7), (2,6), (3,5) \}$ and generic
  multisegment $[1,7]\oplus[2,6]\oplus[3,5]^2\oplus[4,4]$.
\medskip

\centerline{\includegraphics[height=2cm]{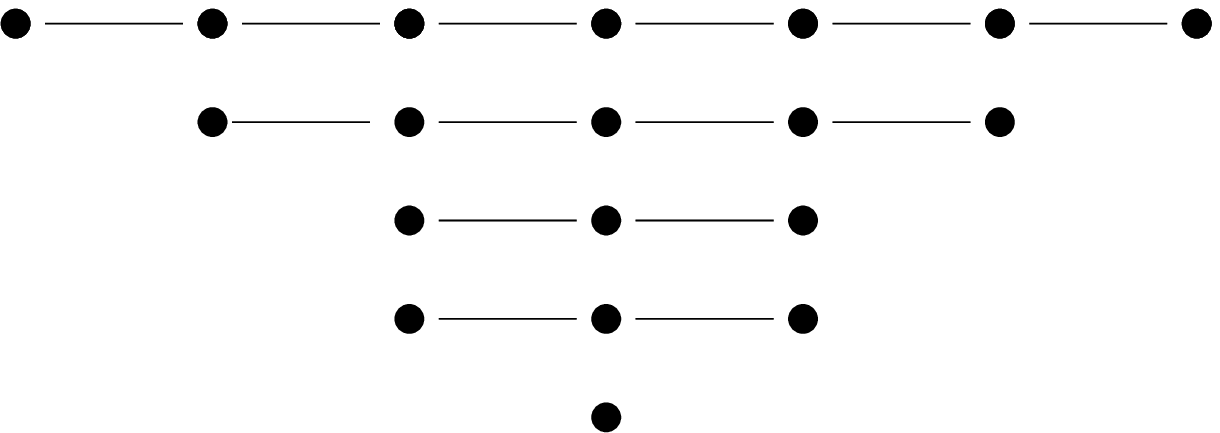}}

\centerline{ {\bf Figure 6.} $M(1,2,4,5,4,2,1)$}
\medskip

and its minimal degenerations look like (all are different and irreducible)

\bigskip

\centerline{\includegraphics[width=10cm]{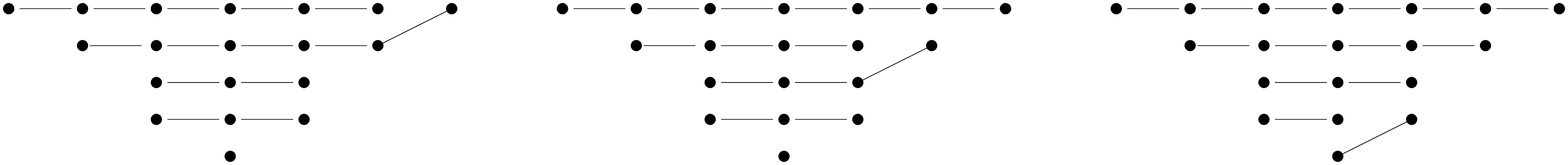}}

\centerline{ {\bf Figure 7.} The minimal degenerations of $M(1,2,4,5,4,2,1)$}
\medskip

\subsection{Pure dimension vectors}

A characterization of all dimension vectors $d$ with $Y$
equidimensional seems to be quite technical. So we restrict the result
to the codimension one case. We define $d$ to be pure, if for alle
elements $(i,j)$ in $I(d)$ we have $d_i = d_j$. The following
proposition is easy to check. 

\begin{Prop}
a) The complement $Y$ of the dense orbit in $\cR(Q,d)$ is equidimensional
of codimension $1$ precisely when $d$ is pure. \\
b) Assume $d$ is pure, then $J(d)$ contains all pairs $(i,i+1)$ and
$I(d)$ only consists of the pairs $(i,j)$ with $d_i = d_j$ and $d_l >
d_i = d_j$ for all $i < l < j$.
\end{Prop}

{\sc Example. }
We have already seen a pure example that is also convex in section
5.3, Example 1. So we
consider a pure one that is not convex. Let $d$ be $(1,2,3,5,3,2,3,2,1)$ 
\medskip

\centerline{\includegraphics[height=2cm]{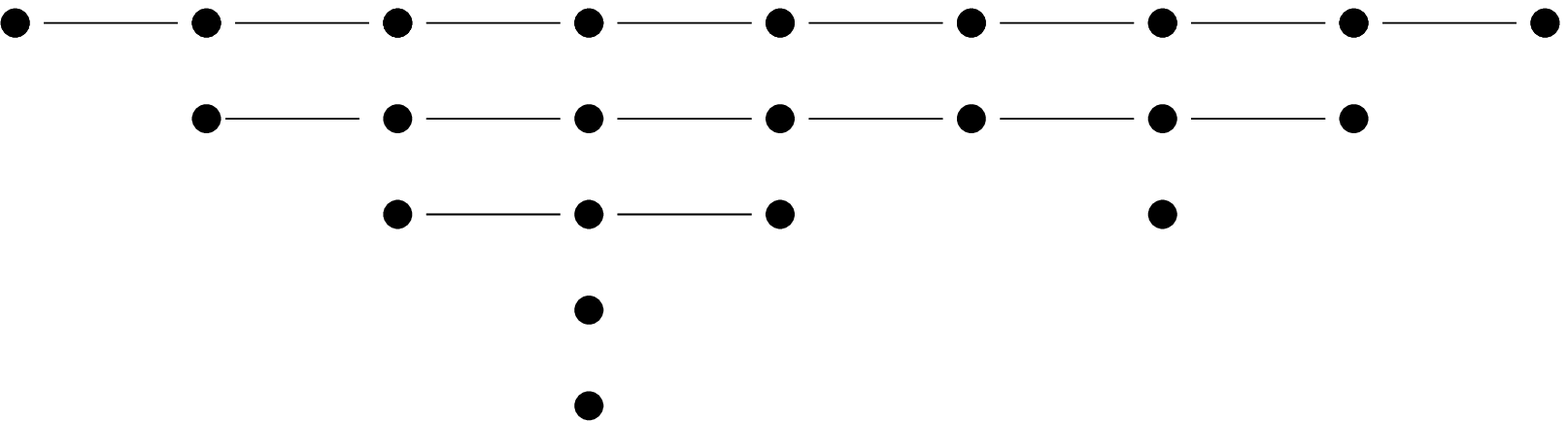}}

\centerline{{\bf Figure 8.} $M(1,2,3,5,3,2,3,2,1)$ }
\medskip

The components are given by $I(1,2,3,5,3,2,3,2,1) =
\{(1,9),(2,6),(6,8),(3,5)\}$. 
\medskip

\centerline{\includegraphics[height=2cm]{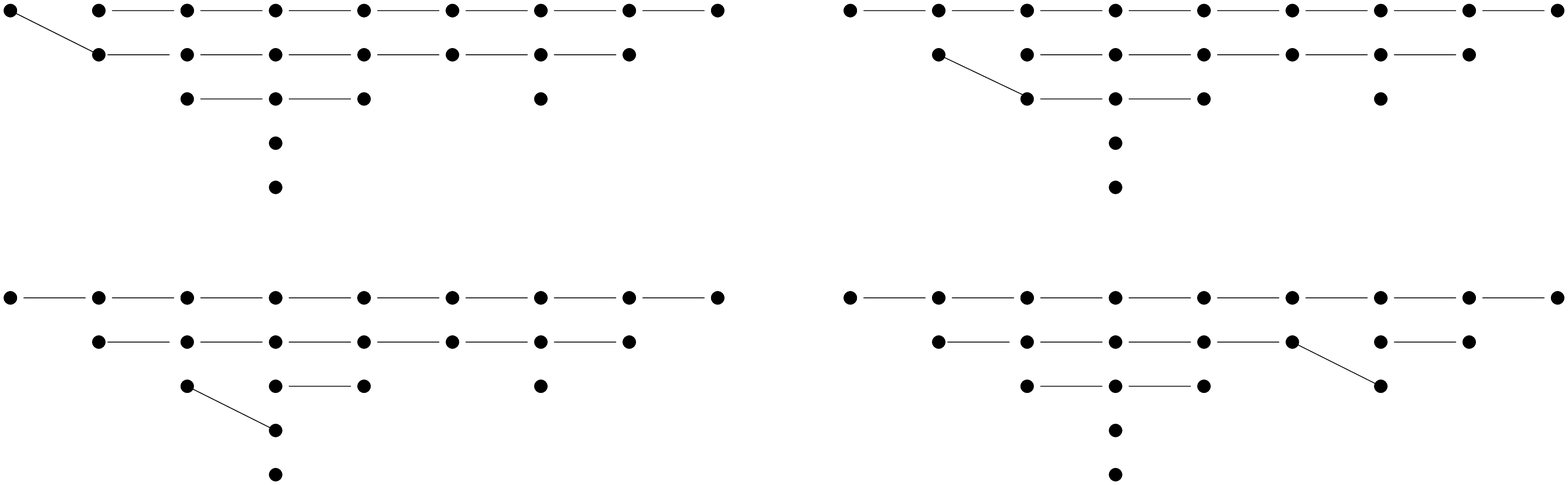}}

\centerline{{\bf Figure 9.} The minimal degenerations of
  $M(1,2,3,5,3,2,3,2,1)$ } 
\medskip





\section{Parabolic group actions}\label{Sparabol}

The results in this note are inspired by the decription of the
complement of the Richardson orbit (the dense orbit) for the action of
a parabolic subgroup in $\mGl_N$ on its unipotent radical as
considered recently in \cite{BH1}. We explain the common idea and some
generalizations. 

\subsection{The Richardson orbit}
Given a dimension vector $d$ as above, then we define a group

$$
P(d) := \{ f \in \mAut(\oplus_{i=1}^t V_i) \mid f(V_j) \subseteq
\oplus_{i=1}^j V_i \mforall j=1,\ldots,t \}
$$

and a vector space

$$
\frp_u(d) := \{ f \in \End(\oplus_{i=1}^t V_i) \mid f(V_j) \subseteq
\oplus_{i=1}^{j - 1} V_i \mforall j=1,\ldots,t \}.
$$

The group $P(d)$ is a standard parabolic subgroup in the General
Linear Group and $\frp_u(d)$ is the Lie algebra of the unipotent radical
of $P(d)$. The group $P(d)$ acts on $\frp_u(d)$, its derived Lie
algebras

$$
\frp_u(d)^{(l)} := \{ f \in \End(\oplus_{i=1}^t V_i) \mid f(V_j) \subseteq
\oplus_{i=1}^{j - l - 1} V_i \mforall j=1,\ldots,t \},
$$ 
and also the quotients $\frp_u(d)/\frp_u(d)^{(l)}$ (and
$\frp_u(d)^{(k)}/\frp_u(d)^{(l)}$ for $k < l$) via conjugation. By a
classical result of Richardson (\cite{Rich}) the group $P(d)$ acts
with a dense orbit on $\frp_u(d)$ and consequently also with a dense
orbit on 
$\frp_u(d)/\frp_u(d)^{(l)}$ for all $l > 1$. In \cite{BH1} we describe
the complement of the dense orbit explicitely using certain rank
conditions on the matrix $A \in \frp_u(d)$. Thus it is desirable to
obtain a similar elementary description of the complement of the dense
orbit in the case $k = l -1$. In fact this case corresponds to a
disjoint union of equioriented Dynkin quivers of type $\Aa$.

The case  $\frp_u(d)/\frp_u(d)^{(l)}$ can be handled using a variation
of the line diagrams introduced in \cite{BH1}.

In contrast, the case $\frp_u(d)^{(l)}$ is still open in general. It
is even not known for general $d$ whether there exists a dense
orbit. A first idea to attack the problem can be found in
\cite{Hvol}. 

\subsection{Irreducible components for the quotients
  $\frp_u(d)/\frp_u(d)^{(l)}$} 

The combinatrorics with line diagrams allows to describe also the
components of the complement of the dense orbit in all quotients
$\frp_u(d)/\frp_u(d)^{(l)}$. For this, we define subvarieties $Z_{i,j}
\subset \frp_u(d)/\frp_u(d)^{(l)}$ by certain rank conditions as in
\cite{BH1}. We claim that there are sets $I^{(l)}(d)
\subseteq J^{(l)}(d)$ so that \\
i) $Z_{i,j}$ is irreducible precisely when $(i,j) \in J^{(l)}(d)$ and
\\
ii) $Z_{i,j}$ is an irreducible component in the complement of the
Richardson orbit, if and only if $(i,j) \in I^{(l)}(d)$. 
\medskip

The definition of $Z_{i,j}$ and the construction of the sets $I^{(l)}(d)
\subseteq J^{(l)}(d)$ can be read off from line diagrams with
connections of length at most $l$. Note that this specializes to the
situation in this note for $l=1$ and to the construction in
\cite{BH1} for $l$ sufficiently large. 

\medskip

{{\small Karin Baur}\\
{\small R\"amistrasse 101}\\
{\small Department of Mathematics}\\
{\small 8092 Z\"urich}\\
{\small Switzerland}\\
{\small E-mail: baur@math.ethz.ch}\\
{\small http://www.math.ethz.ch/$\sim$baur}}
\medskip

{{\small Lutz Hille}\\
{\small Mathematisches Institut}\\
{\small Universit\"at M\"unster}\\
{\small Einsteinstr.~62}\\
{\small 48149 M\"unster}\\
{\small Germany}\\
{\small E-mail: lutz.hille@uni-muenster.de}\\
{\small http://wwwmath.uni-muenster.de/reine/u/lutz.hille/}}


\begin{thebibliography}{100000}

\bibitem[AF]{AbeasisdelFra}
Abeasis, S.; Del Fra, A. Degenerations for the representations of an
equioriented quiver of type $A\sb{m}$.  Boll. Un. Mat. Ital. Suppl.
1980,  no. 2, 157--171.

\bibitem[AFK]{AbeasisdelFraKraft} Abeasis, S.; Del Fra, A.; Kraft,
  H. {\em The geometry of representations of $A\sb{m}$. }
  Math.~Ann.~{\bf   256}
  (1981), no.~3, 401--418. 

\bibitem[BH]{BH1} Baur, K.; Hille, L. {\em On the Complement
  of Richardson Orbit}, preprint 2009

\bibitem[BHRR]{BHRR} Br\"ustle, T.; Hille, L.; Ringel, C.~M.; R\"ohrle, G.
{\em The $\Delta$-filtered modules without self-extensions for the
  Auslander algebra of $k[T]/\langle T\sp n\rangle$. }
Algebr.~Represent.~Theory {\bf  2}  (1999),  no.~3, 295--312.

\bibitem[BHRZ]{BHRZ}
Br\"ustle, T.; Hille, L.; R\"ohrle, G.; Zwara, G. {\em The
  Bruhat-Chevalley order of parabolic group actions in general linear
  groups and degeneration for $\Delta$-filtered modules. }  Adv.~Math.~{\bf
  148}  (1999),  no.~2, 203--242.
 
\bibitem[F]{Fultontoric} Fulton, W. {\em Introduction to toric
    varieties. }  Annals of Mathematics Studies, {\bf 131}. The William
    H.~Roever Lectures in Geometry. Princeton University Press,
    Princeton 
 
\bibitem[H1]{Hvol} Hille, L. {\em  On the volume of a tilting
  module. }  Abh.~Math.~Sem.~Univ.~Hamburg {\bf 76}  (2006), 261--277.

\bibitem[H2]{Hhabil} Hille, L. {\em Aktionen algebraischer
  Gruppen, geometrische Quotienten und K\"ocher. } habilitation
  thesis, Hamburg 2003.


\bibitem[KZ]{KnightZelevinsky} 
Knight, H.; Zelevinsky, A.
{\em Representations of quivers of type $A$ and the multisegment
  duality. }
Adv.~Math.~{\bf 117} (1996), no.~2, 273--293

\bibitem[P]{Piasetsky}
Pjaseckii, V. S. {\em Linear Lie groups that act with a finite number
  of orbits. } Funkcional.~Anal.~i Prilozen.~{\bf 9}  (1975), no.~4,
85--86. 

\bibitem[Rc]{Rich} Richardson, R. W., Jr. {\em  Conjugacy classes in
  parabolic subgroups of semisimple algebraic groups. }  Bull.~London
  Math.~Soc.~{\bf  6}  (1974), 21--24.


\bibitem[RS]{RiedtmannSchofield2}
Riedtmann, C., Schofield, A. {\em
On a simplicial complex associated with tilting modules. }
Comment.~Math.~Helv.~{\bf 66} (1991), no.~1, 70--78.

\bibitem[Rn]{Ringel}
Ringel, C.~M. {\em  The multisegment duality and the preprojective
  algebras of type $A$. } AMA Algebra Montp.~Announc.~1999 (electr.).

\bibitem[S]{Spaltenstein} Spaltenstein, N.
{\em Classes unipotentes et sous-groupes de Borel. }
Lecture Notes in Mathematics, {\bf 946}. Springer-Verlag, Berlin-New
York, 1982.
 
\bibitem[Z]{Zwara}  Zwara, G. {\em Degenerations for modules over
    representation-finite algebras. }  Proc.~Amer.~Math.~Soc.~{\bf 127}  no.~5
    (1999).
\end{thebibliography}
\end{document}